\documentclass[10pt]{amsart}

\usepackage{amsmath,amssymb}
\usepackage{url}

\newtheorem{theorem}{Theorem}[section]
\newtheorem{lemma}[theorem]{Lemma}
\newtheorem{corollary}[theorem]{Corollary}
\theoremstyle{remark}
\newtheorem{remark}[theorem]{Remark}
\theoremstyle{definition}
\newtheorem{definition}[theorem]{Definition}

\newcommand{\ped}{\operatorname{ped}}
\newcommand{\op}{\operatorname{op}}

\title[Conjecture of Banerjee and Bringmann]
{A solution to a conjecture of Banerjee and Bringmann on restricted overpartitions}
\author{Mohamed El Bachraoui}
\address{Department of Mathematical Sciences, UAE University, Al Ain, United Arab Emirates}
\email{melbachraoui@uaeu.ac.ae}
\date{}

\begin{document}
\begin{abstract}
In this paper, we confirm a conjecture of Banerjee and Bringmann on the parity
of a restricted overpartition function. Our main result is the congruence
\[
  \op^{[1]}_1(4n)\equiv\ped(n)\pmod2,
\]
where $\ped(n)$ counts partitions of $n$ in which the even parts are
distinct. This gives a triangular-number characterization of the parity
of the overpartition function and proves both directions of the
conjectured Legendre-symbol criterion. We also record an extension to
the corresponding family of restricted overpartition functions.
\end{abstract}
\maketitle

\section{Introduction and main results}

Banerjee and Bringmann \cite{BB} introduced restricted partition and
overpartition functions whose generating functions are connected with
universal mock theta functions. They also conjectured a Legendre-symbol
criterion for the parity of one of their overpartition functions. We
settle this conjecture by relating its coefficients, modulo $2$, to
partitions with distinct even parts. Our proof uses elementary identities
for $q$-series.

An overpartition is a partition in which the first occurrence of a part
may be overlined; see Corteel and Lovejoy \cite{CL}. We use the notation
\[
 (a;q)_m=\prod_{j=0}^{m-1}(1-aq^j),
\]
with $(a;q)_0=1$ and the usual interpretation when $m=\infty$.

\begin{definition}\label{def:op1}
Following Banerjee and Bringmann \cite{BB}, let $\op^{[1]}_1(n)$ count
the overpartitions of $n$ satisfying the following restrictions:
\begin{enumerate}
\item every even part is non-overlined;
\item the largest part is even;
\item every odd part occurs at most twice, and if it occurs twice, one
      occurrence is overlined.
\end{enumerate}
For $n=0$, we set $\op^{[1]}_1(0)=1$, corresponding to the empty
overpartition. The generating function obtained in
\cite[Proof of Theorem 1.2]{BB} is
\begin{equation}\label{eq:BB-gf}
 \sum_{n\geq0}\op^{[1]}_1(n)q^n
 =\sum_{m\geq0}
 \frac{q^{2m}(-q;q^2)_m^2}{(q^2;q^2)_m}.
\end{equation}
\end{definition}

In \cite[Conjecture 4.1]{BB}, Banerjee and Bringmann
proposed the following criterion: for every odd prime $p$ and every
integer $k\geq0$,
\begin{equation}\label{eq:original-conjecture}
  \op^{[1]}_1(4pn+4k)\equiv0\pmod2\quad\text{for every }n\geq0
  \quad\Longleftrightarrow\quad
  \left(\frac{8k+1}{p}\right)=-1.
\end{equation}
Here $\left(\frac{\cdot}{p}\right)$ denotes the Legendre symbol.

To state our result, let $\ped(n)$ denote the number of partitions of $n$
in which even parts are distinct and odd parts are unrestricted. The generating function for $\ped(n)$ is
\begin{equation}\label{eq:ped-gf}
 \sum_{n\geq0}\ped(n)q^n
 =\frac{(-q^2;q^2)_\infty}{(q;q^2)_\infty}
 =\frac{(q^4;q^4)_\infty}{(q;q)_\infty}.
\end{equation}
Thus $\ped(n)$ also counts the $4$-regular partitions of $n$, that is,
partitions with no part divisible by $4$.

The arithmetic of $\ped(n)$ has been studied extensively; see, for
example, Andrews, Hirschhorn, and Sellers \cite{AHS}. In particular,
the parity of $\ped(n)$ is governed by triangular numbers, as follows
from Gauss's triangular-number identity. Bessenrodt and Pak
\cite[Section~6]{BP} also noted this
parity consequence and studied related partition classes by involutions.
Andrews and El Bachraoui \cite{AE} subsequently investigated identities
relating restricted partitions with distinct even parts and restricted
$4$-regular partitions. In another paper, Andrews and El Bachraoui
\cite{AE2} obtained a modulo-$4$ description for a two-color partition
function whose coefficients also have triangular-number parity, together
with a related criterion for the vanishing of arithmetic progressions.

The new connection here is between $\ped(n)$ and the restricted
overpartition function of Banerjee and Bringmann. Our main result is the
following coefficientwise congruence.
\begin{theorem}\label{thm:main}
 For every integer $n\geq0$,
 \begin{equation}\label{eq:main}
   \op^{[1]}_1(4n)\equiv\ped(n)\pmod2.
 \end{equation}
 In particular,
 \begin{equation}\label{eq:triangular}
 \op^{[1]}_1(4n)\equiv
 \begin{cases}
  1,&\text{if }n=j(j+1)/2\text{ for some }j\geq0,\\
  0,&\text{otherwise}
 \end{cases}\pmod2
 \end{equation}
\end{theorem}

The congruence \eqref{eq:main} is the main assertion; the
triangular-number characterization follows from the classical parity of
$\ped(n)$. Together they prove Conjecture 4.1 in both directions.
\begin{corollary}\label{cor:conjecture}
 Let $p$ be an odd prime and $k$ a nonnegative integer. Then
 \[
  \op^{[1]}_1(4pn+4k)\equiv0\pmod2\quad\text{for all }n\geq0
  \quad\Longleftrightarrow\quad
  \left(\frac{8k+1}{p}\right)=-1.
 \]
\end{corollary}

We also record the extension to the statistic $\op^{[1]}_r$ of \cite{BB}.
By its definition, all parts are divisible by $r$, and division by $r$
identifies these overpartitions with those counted by $\op^{[1]}_1$.
\begin{corollary}\label{cor:r}
 For all positive integers $r$ and all $n\geq0$,
 \[
  \op^{[1]}_r(4rn)\equiv\ped(n)\pmod2.
 \]
 In particular, for every odd prime $p$ and $k\geq0$,
 \[
  \op^{[1]}_r(4rpn+4rk)\equiv0\pmod2\quad\text{for all }n\geq0
  \quad\Longleftrightarrow\quad
  \left(\frac{8k+1}{p}\right)=-1.
 \]
\end{corollary}

\section{Proof of Theorem \ref{thm:main}}

The main step is a parity identity for the complete generating function
of $\op^{[1]}_1(n)$. 

\begin{lemma}\label{lem:bilateral}
We have
\begin{equation}\label{eq:bilateral}
 \sum_{n\geq0}\op^{[1]}_1(n)q^n
 \equiv\frac{1}{(q^4;q^4)_\infty}
 \sum_{j\in\mathbb Z}q^{6j^2+4j}\pmod2.
\end{equation}
\end{lemma}

\begin{proof}
We start from \eqref{eq:BB-gf}. All congruences in this proof are
identities of formal power series with coefficients reduced modulo $2$.
Since
\[
 (-q;q^2)_m^2
 =\prod_{j=0}^{m-1}(1+q^{2j+1})^2
 \equiv\prod_{j=0}^{m-1}(1+q^{4j+2})
 =(-q^2;q^4)_m\pmod2,
\]
we obtain
\begin{equation}\label{eq:first-reduction}
 \sum_{n\geq0}\op^{[1]}_1(n)q^n
 \equiv\sum_{m\geq0}
 \frac{q^{2m}(-q^2;q^4)_m}{(q^2;q^2)_m}
 \pmod2.
\end{equation}

We next separate the terms with even and odd $m$. The denominator
factors as
\begin{equation}\label{eq:denominator-factors}
 \begin{aligned}
 (q^2;q^2)_{2l}&=(q^2;q^4)_l(q^4;q^4)_l,\\
 (q^2;q^2)_{2l+1}&=(q^2;q^4)_{l+1}(q^4;q^4)_l.
 \end{aligned}
\end{equation}
Splitting the numerator products after their first $l$ or $l+1$
factors, respectively, gives
\begin{equation}\label{eq:numerator-factors}
 \begin{aligned}
 (-q^2;q^4)_{2l}
   &=(-q^2;q^4)_l(-q^{4l+2};q^4)_l,\\
 (-q^2;q^4)_{2l+1}
   &=(-q^2;q^4)_{l+1}(-q^{4l+6};q^4)_l.
 \end{aligned}
\end{equation}
Since $(-q^2;q^4)_j\equiv(q^2;q^4)_j\pmod2$ for all $j\geq0$,
we may cancel the corresponding factors in
\eqref{eq:denominator-factors} and \eqref{eq:numerator-factors}.
For even indices this yields
\begin{equation}\label{eq:even-m}
 \frac{(-q^2;q^4)_{2l}}{(q^2;q^2)_{2l}}
 \equiv
 \frac{(-q^{4l+2};q^4)_l}{(q^4;q^4)_l}\pmod2,
\end{equation}
while for odd indices it gives
\begin{equation}\label{eq:odd-m}
 \frac{(-q^2;q^4)_{2l+1}}{(q^2;q^2)_{2l+1}}
 \equiv
 \frac{(-q^{4l+6};q^4)_l}{(q^4;q^4)_l}\pmod2.
\end{equation}
Every cancelled factor has constant term $1$, so cancellation is
legitimate in the formal power series ring modulo $2$.
Substituting \eqref{eq:even-m} and \eqref{eq:odd-m} into
\eqref{eq:first-reduction}, we find
\begin{equation}\label{eq:parity-m}
 \begin{aligned}
 \sum_{n\geq0}\op^{[1]}_1(n)q^n
 \equiv{}&\sum_{l\geq0}
 \frac{q^{4l}(-q^{4l+2};q^4)_l}{(q^4;q^4)_l}\\
 &+\sum_{l\geq0}
 \frac{q^{4l+2}(-q^{4l+6};q^4)_l}{(q^4;q^4)_l}
 \pmod2.
 \end{aligned}
\end{equation}

Recall that the Gaussian binomial coefficient is given by
\begin{equation}\label{eq:gaussian-binomial}
 \begin{bmatrix}l\\j\end{bmatrix}_q
 =\frac{(q;q)_l}{(q;q)_j(q;q)_{l-j}}
 \qquad (0\leq j\leq l).
\end{equation}
The finite $q$-binomial theorem \cite[(3.3.6)]{AndrewsBook},
with nome $q^4$, gives
\begin{equation}\label{eq:finite-q-binomial}
 (-q^a;q^4)_l
 =\sum_{j=0}^{l}q^{aj+2j(j-1)}
 \begin{bmatrix}l\\j\end{bmatrix}_{q^4}
 \qquad (a\geq1).
\end{equation}
We apply \eqref{eq:finite-q-binomial} with $a=4l+2$ and
$a=4l+6$. Including the leading powers of $q$ in
\eqref{eq:parity-m}, the resulting exponents are
\begin{equation}\label{eq:exponents}
 4l+4lj+2j^2,\qquad 4l+2+4lj+4j+2j^2
 \qquad(0\leq j\leq l).
\end{equation}
Thus \eqref{eq:parity-m}, \eqref{eq:gaussian-binomial}, and
\eqref{eq:finite-q-binomial} give
\begin{equation}\label{eq:lemma-double-sums}
 \begin{aligned}
 \sum_{n\geq0}\op^{[1]}_1(n)q^n
 \equiv{}&\sum_{l\geq0}\sum_{j=0}^l
 \frac{q^{4l+4lj+2j^2}}
 {(q^4;q^4)_j(q^4;q^4)_{l-j}}\\
 &+\sum_{l\geq0}\sum_{j=0}^l
 \frac{q^{4l+2+4lj+4j+2j^2}}
 {(q^4;q^4)_j(q^4;q^4)_{l-j}}
 \pmod2.
 \end{aligned}
\end{equation}

Instead of extracting the powers of $q^4$ at this point, put
$l=j+m$ in both double sums. The exponents in
\eqref{eq:exponents} become
\[
 \begin{aligned}
 4l+4lj+2j^2&=6j^2+4j+4(j+1)m,\\
 4l+2+4lj+4j+2j^2&=6j^2+8j+2+4(j+1)m.
 \end{aligned}
\]
Consequently, \eqref{eq:lemma-double-sums} becomes
\begin{equation}\label{eq:lemma-iterated-sums}
 \begin{aligned}
 \sum_{n\geq0}\op^{[1]}_1(n)q^n
 \equiv{}&\sum_{j\geq0}
 \frac{q^{6j^2+4j}}{(q^4;q^4)_j}
 \sum_{m\geq0}\frac{q^{4(j+1)m}}{(q^4;q^4)_m}\\
 &+\sum_{j\geq0}
 \frac{q^{6j^2+8j+2}}{(q^4;q^4)_j}
 \sum_{m\geq0}\frac{q^{4(j+1)m}}{(q^4;q^4)_m}
 \pmod2.
 \end{aligned}
\end{equation}

Euler's special case of the $q$-binomial theorem
\cite[(2.2.5)]{AndrewsBook} states that
\begin{equation}\label{eq:euler-sum}
 \sum_{m\geq0}\frac{q^{4am}}{(q^4;q^4)_m}
 =\frac{1}{(q^{4a};q^4)_\infty}
 \qquad(a\geq1).
\end{equation}
We also have the product identity
\begin{equation}\label{eq:lemma-product}
 (q^4;q^4)_j(q^{4j+4};q^4)_\infty
 =(q^4;q^4)_\infty.
\end{equation}
Applying \eqref{eq:euler-sum} with $a=j+1$ to the inner sums in
\eqref{eq:lemma-iterated-sums}, and then using
\eqref{eq:lemma-product}, gives
\begin{equation}\label{eq:lemma-theta-half}
 \begin{aligned}
 \sum_{n\geq0}\op^{[1]}_1(n)q^n
 &\equiv\sum_{j\geq0}
 \frac{q^{6j^2+4j}}
 {(q^4;q^4)_j(q^{4j+4};q^4)_\infty}\\
 &\quad+\sum_{j\geq0}
 \frac{q^{6j^2+8j+2}}
 {(q^4;q^4)_j(q^{4j+4};q^4)_\infty}\\
 &=\frac{1}{(q^4;q^4)_\infty}
 \left(\sum_{j\geq0}q^{6j^2+4j}
 +\sum_{j\geq0}q^{6j^2+8j+2}\right)
 \pmod2.
 \end{aligned}
\end{equation}
Finally,
\[
 6j^2+8j+2=6(-j-1)^2+4(-j-1).
\]
The first sum in \eqref{eq:lemma-theta-half} corresponds to
nonnegative indices and the second to negative indices. These ranges
are disjoint and together cover $\mathbb Z$, proving
\eqref{eq:bilateral}.
\end{proof}

\begin{proof}[Proof of Theorem \ref{thm:main}]
The denominator in \eqref{eq:bilateral} involves only powers of
$q^4$. Moreover,
\[
 6j^2+4j\equiv2j^2\pmod4,
\]
so $6j^2+4j$ is divisible by $4$ if and only if $j$ is even.
Extracting the terms of degree divisible by $4$ in
\eqref{eq:bilateral} and writing $j=2k$, we obtain
\begin{equation}\label{eq:even-extraction}
 \sum_{n\geq0}\op^{[1]}_1(4n)q^{4n}
 \equiv\frac{1}{(q^4;q^4)_\infty}
 \sum_{k\in\mathbb Z}q^{24k^2+8k}\pmod2.
\end{equation}
Euler's pentagonal number theorem \cite[(1.3.1)]{AndrewsBook},
applied with nome $q^{16}$ and with the bilateral summation
index replaced by its negative, gives
\begin{equation}\label{eq:pentagonal-identity}
 \begin{aligned}
 (q^{16};q^{16})_\infty
 &=\sum_{k\in\mathbb Z}(-1)^k q^{8k(3k+1)}\\
 &\equiv\sum_{k\in\mathbb Z}q^{24k^2+8k}\pmod2.
 \end{aligned}
\end{equation}
Combining \eqref{eq:even-extraction} and
\eqref{eq:pentagonal-identity}, we obtain
\begin{equation}\label{eq:ped-extraction}
 \sum_{n\geq0}\op^{[1]}_1(4n)q^{4n}
 \equiv\frac{(q^{16};q^{16})_\infty}
 {(q^4;q^4)_\infty}
 =\sum_{n\geq0}\ped(n)q^{4n}\pmod2.
\end{equation}
The final equality follows from \eqref{eq:ped-gf}.
Comparing coefficients of $q^{4n}$ proves \eqref{eq:main}.

To obtain the triangular-number description, note that
\[
 (q^2;q^2)_\infty^2\equiv(q^4;q^4)_\infty\pmod2
\]
and recall Gauss's triangular-number identity \cite[(2.2.13)]{AndrewsBook} stating that
\begin{equation}\label{eq:gauss-identity}
 \frac{(q^2;q^2)_\infty^2}{(q;q)_\infty}
 =\sum_{j\geq0}q^{j(j+1)/2}.
\end{equation}
Consequently, \eqref{eq:ped-gf} and \eqref{eq:gauss-identity}
give
\begin{equation}\label{eq:gauss}
 \begin{aligned}
 \sum_{n\geq0}\ped(n)q^n
 &\equiv\frac{(q^2;q^2)_\infty^2}{(q;q)_\infty}\\
 &=\sum_{j\geq0}q^{j(j+1)/2}\pmod2.
 \end{aligned}
\end{equation}
Together with \eqref{eq:main}, this proves \eqref{eq:triangular}
and completes the proof of Theorem \ref{thm:main}.
\end{proof}

\section{Proof of the corollaries}

\begin{proof}[Proof of Corollary \ref{cor:conjecture}]
By Theorem \ref{thm:main}, $\op^{[1]}_1(4pn+4k)$ is odd exactly when
$pn+k$ is triangular, that is, when
\begin{equation}\label{eq:square}
 8pn+8k+1=(2j+1)^2
\end{equation}
for some integers $n,j\geq0$.

Suppose first that
\[
 \left(\frac{8k+1}{p}\right)=-1.
\]
Reducing \eqref{eq:square} modulo $p$ would give
$(2j+1)^2\equiv8k+1\pmod p$, which is impossible.
Hence $\op^{[1]}_1(4pn+4k)$ is even for every $n\geq0$.

Conversely, suppose that the Legendre symbol is $0$ or $1$.
There exists an integer $a$ satisfying
\[
 a^2\equiv8k+1\pmod p.
\]
Since $p$ is odd, the integers congruent to $a$ modulo $p$ include
arbitrarily large positive odd integers: adding $p$ changes the parity.
Choose one such integer $m$. Then
\[
 m^2\equiv8k+1\pmod p,
 \qquad m^2\equiv1\equiv8k+1\pmod8.
\]
Since $8$ and $p$ are coprime, $8p$ divides
$m^2-(8k+1)$. Choosing $m$ sufficiently large also ensures that
\[
 n:=\frac{m^2-(8k+1)}{8p}\geq0.
\]
Writing $m=2j+1$ with $j\geq0$ gives \eqref{eq:square}.
This construction also applies when the Legendre symbol is $0$, since
we may then take $a=0$.
Thus $pn+k=j(j+1)/2$, so Theorem \ref{thm:main} yields
$\op^{[1]}_1(4pn+4k)\equiv1\pmod2$ for this value of $n$.
Therefore the coefficients in the progression cannot all be even.
This proves the converse and the corollary.
\end{proof}

\begin{proof}[Proof of Corollary \ref{cor:r}]
Every part of an overpartition counted by $\op^{[1]}_r(rn)$ is a
multiple of $r$. Dividing all its parts by $r$ sends parts divisible by
$2r$ to even parts and parts congruent to $r$ modulo $2r$ to odd parts.
It preserves which occurrences are overlined, the multiplicity
restrictions, and the requirement that the largest part be divisible by
$2r$. Conversely, multiplication of every part by $r$ reverses this
operation. Hence this is a bijection with the overpartitions counted by
$\op^{[1]}_1(n)$.
Consequently,
\[
 \op^{[1]}_r(rn)=\op^{[1]}_1(n)\qquad(n\geq0).
\]
The assertions now follow from Theorem \ref{thm:main} and
Corollary \ref{cor:conjecture}.
\end{proof}

\begin{remark}
A direct combinatorial proof of \eqref{eq:main} would be of interest.
\end{remark}

\enlargethispage{4\baselineskip}

\end{document}